# RECURSIVE ESTIMATION OF TIME-AVERAGE VARIANCE CONSTANTS[1]

By Wei Biao Wu

*University of Chicago*

For statistical inference of means of stationary processes, one needs to estimate their time-average variance constants (TAVC) or long-run variances. For a stationary process, its TAVC is the sum of all its covariances and it is a multiple of the spectral density at zero. The classical TAVC estimate which is based on batched means does not allow recursive updates and the required memory complexity is $O(n)$. We propose a faster algorithm which recursively computes the TAVC, thus having memory complexity of order $O(1)$ and the computational complexity scales linearly in $n$. Under short-range dependence conditions, we establish moment and almost sure convergence of the recursive TAVC estimate. Convergence rates are also obtained.

**1. Introduction.** Let $(X_i)_{i\in\mathbb{Z}}$ be a stationary and ergodic process with mean $\mu = \mathbb{E}(X_0)$ and finite variance; let $\gamma(k) = \text{cov}(X_0, X_k)$, $k \in \mathbb{Z}$, be the covariance function. Given the observations $X_1, \ldots, X_n$, a simple estimate of $\mu$ is the sample mean $\bar{X}_n = n^{-1} \sum_{i=1}^{n} X_i$. Under suitable conditions on $(X_i)$, $\bar{X}_n$ is asymptotically normal:

$$(1) \qquad n^{1/2}(\bar{X}_n - \mu) = n^{-1/2} \sum_{i=1}^{n} (X_i - \mu) \Rightarrow N(0, \sigma^2),$$

where $\Rightarrow$ denotes convergence in distribution and $\sigma^2$ is called the time-average variance constant (TAVC), long-run variance or asymptotic variance parameter. Goodman and Sokal (1989) called $\sigma^2/\gamma(0)$ the integrated autocorrelation time. There exists a huge literature on the central limit theory for stationary processes. See, for example, Ibragimov and Linnik (1971) and Bradley (2007).

Received June 2007; revised November 2008.

[1]Supported in part by NSF Grant DMS-04-78704.

To conduct statistical inference for $\mu$, one needs to estimate $\sigma^2$. Under suitable conditions, $\sigma^2 = \sum_{k \in \mathbb{Z}} \gamma(k)$. The estimation of $\sigma^2$ is an important problem in statistical inference of time series and it has a long history. Given $X_1, \ldots, X_n$, let $1 \leq l_n \leq n$ be the block length satisfying $l_n \to \infty$ and $l_n/n \to 0$. Based on the batched means $\sum_{i=j}^{j+l_n-1} X_i/l_n$, $1 \leq j \leq n - l_n + 1$, one can estimate $\sigma^2$ by

$$(2) \qquad \sigma_n^2(l_n) = \frac{l_n}{n - l_n + 1} \sum_{j=1}^{n-l_n+1} \left( \frac{1}{l_n} \sum_{i=j}^{j+l_n-1} X_i - \bar{X}_n \right)^2.$$

The estimate $\sigma_n^2(l_n)$ appears in several contexts and it is closely related to Bartlett's spectral density estimate. As an alternative, one can propose a similar estimate by using the nonoverlapped batched means $\sum_{i=j}^{j+l_n-1} X_i/l_n$, $j = 1, 1 + l_n, 1 + 2l_n, \ldots$. Asymptotic properties of $\sigma_n^2(l_n)$ have been extensively studied; see, for example, Alexopoulos and Goldsman (2004), Song and Schmeiser (1995), Bühlmann (2002), Lahiri (2003), Politis, Romano and Wolf (1999) and Jones et al. (2006), among others. For other works on estimation of $\sigma^2$, Chauveau and Diebolt (2003) used multiple parallel chains, and Robert (1995) considered Harris recurrent chains. The estimation of $\sigma^2$ is related to the problem of Markov chain Monte Carlo (MCMC) convergence assessment; see Brooks and Roberts (1998), Chauveau and Diebolt (1999) and Chauveau, Diebolt and Robert (1998), among others.

It is well known that $\bar{X}_n$ can be recursively computed in the sense that, if a new observation $X_{n+1}$ is available, then $\bar{X}_{n+1}$ can be computed as $(n\bar{X}_n + X_{n+1})/(n+1)$. Hence, the memory complexity for computing $\bar{X}_n$ is $O(1)$. However, this nice property is no longer present in the estimate $\sigma_n^2(l_n)$ in (2). There is no simple algebraic relation between $\sigma_{n+1}^2(l_{n+1})$ and $\sigma_n^2(l_n)$. To compute $\sigma_{n+1}^2(l_{n+1})$, if $l_n \neq l_{n+1}$, one then has to update all batched means and the memory complexity is $O(n)$. In computationally intensive problems, it is desirable to have a recursive estimate. For example, in MCMC experiments, one sequentially generates $X_1, X_2, \ldots$. At each stage, based on (1), a $(1 - \alpha)$ confidence interval of $\mu$ can be constructed as $\bar{X}_n \pm z_{1-\alpha/2}\hat{\sigma}_n/\sqrt{n}$, where $z_{1-\alpha/2}$ is the $(1-\alpha/2)$th percentile of a standard normal distribution, $0 < \alpha < 1$. As argued in Geyer (1992), Fishman (1996) and Jones et al. (2006), among others, for convergence diagnostics of Markov chain Monte Carlo algorithms, one can terminate the simulation by choosing $n$ such that the interval is sufficiently small. Quick update of $\hat{\sigma}_n$ is essential for efficient sequential monitoring and testing. For example, to test the hypothesis $\mu = \mu_0$, we can consider the test statistic $\sqrt{n}|\bar{X}_n - \mu_0|/\hat{\sigma}_n$, which can be quickly calculated via sequentially updating.

A common practice in MCMC simulations is to run multiple i.i.d. copies of the chain. One can run, for example, 100 copies of the chain and then conduct convergence diagnostics based on comparison of asymptotic variances of each chain. In such cases the computational and memory advantage of our recursive algorithm is more appealing.

The rest of the paper is structured as follows. A sequential estimate $\hat{\sigma}_n^2$ of $\sigma^2$ is introduced in Section 2. Namely, at each stage $n$, $\hat{\sigma}_n^2$ can be updated within $O(1)$ steps so that the computational complexity scales linearly in $n$. The moment and almost sure convergence properties are presented in Section 3 and some implementation issues are discussed in Section 4. Section 5 provides applications to Markov chains and linear processes. Proofs are given in the Appendix.

We now introduce some notation. A random variable $\xi$ is said to be in $\mathcal{L}^p$ ($p > 0$) if $\|\xi\|_p := [\mathbb{E}(|\xi|^p)]^{1/p} < \infty$. Write $\|\xi\| = \|\xi\|_2$. For two real sequences $(a_n)$ and $(b_n)$, write $a_n \sim b_n$ if $\lim_{n\to\infty} a_n/b_n = 1$ and $a_n \asymp b_n$ if there exists a constant $c > 0$ such that $1/c \leq |a_n/b_n| \leq c$ for all large $n$. Let $S_n = X_1 + \cdots + X_n - n\mu$ and $S_n^* = \max_{i \leq n} |S_i|$.

**2. Recursive TAVC estimates.** For ease of reading, we assume at the outset that $\mu = 0$. To define our recursive TAVC estimate, let $(a_k)_{k\in\mathbb{N}}$ be a strictly increasing integer-valued sequence such that $a_1 = 1$ and $a_{k+1} - a_k \to \infty$ as $k \to \infty$. Based on $(a_k)_{k\in\mathbb{N}}$, define another sequence $(t_i)_{i\in\mathbb{N}}$ as $t_i = a_k$ if $a_k \leq i < a_{k+1}$. As a simple example, let $a_k = k^2$. Then $t_i = \lfloor \sqrt{i} \rfloor^2$, where $\lfloor u \rfloor = \max\{k \in \mathbb{Z} : k \leq u\}$ is the integer part of $u$. Given $X_1, \ldots, X_n$, define

$$(3) \qquad V_n = \sum_{i=1}^n W_i^2 \qquad \text{where } W_i = X_{t_i} + X_{t_i+1} + \cdots + X_i,$$

and

$$(4) \qquad v_n = \sum_{i=1}^n l_i \qquad \text{where } l_i = i - t_i + 1.$$

We propose to estimate the TAVC $\sigma^2$ by $V_n/v_n$. In the estimate (2), for a given $n$, the block size $l_n$ is the same for different blocks. In $V_n$, however, different blocks have different block lengths. Let $B_k = \{a_k, a_k + 1, \ldots, a_{k+1} - 1\}$. Assume $a_k \leq n < a_{k+1}$. Then $t_n = a_k$ and $W_n = X_{a_k} + X_{a_k+1} + \cdots + X_n$. If $n+1 \neq t_{n+1}$, then $n+1$ still belongs to the block $B_k$ and $W_{n+1} = W_n + X_{n+1}$. On the other hand, however, if $n + 1 = t_{n+1}$, then $t_{n+1} = a_{k+1}$ and $n + 1$ belongs to the next block $B_{1+k}$, and $W_{n+1}$ now becomes $X_{n+1}$. Combining these two cases, we have $W_{n+1} = W_n \mathbf{1}_{n+1 \neq t_{n+1}} + X_{n+1}$. For $n \in \mathbb{N}$, choose $k_n \in \mathbb{N}$ such that $a_{k_n} \leq n < a_{1+k_n}$. Then $a_{k_n} = t_n$. To summarize, we propose the following recursive algorithm:

ALGORITHM 1. At stage $n$, we store $(n, k_n, a_{k_n}, v_n, V_n, W_n)$. Note that $t_n = a_{k_n}$. At stage $n + 1$, we update the vector by:

1. If $n+1 = a_{1+k_n}$, let $k_{n+1} = 1 + k_n$ and $W_{n+1} = X_{n+1}$; If $n+1 \neq a_{1+k_n}$, then let $k_{n+1} = k_n$ and $W_{n+1} = W_n + X_{n+1}$,
2. $V_{n+1} = V_n + W_{n+1}^2$,
3. $v_{n+1} = v_n + (n+2 - t_{n+1})$, where $t_{n+1} = a_{k_{n+1}}$.

Output: $\hat{\sigma}_{n+1}^2 = V_{n+1}/v_{n+1}$.

To implement Algorithm 1, one needs to specify the sequence $(a_k)_{k\geq 1}$. A simple choice is that $a_k = \lfloor ck^p \rfloor$, $k \geq 1$, where $c > 0$ and $p > 1$ are constants (cf. Remark 2 and Theorem 2). We now compute $t_n$ for the sequence $(a_k)_{k\geq 1}$. To this end, let $k \in \mathbb{N}$ be such that $a_k \leq n < a_{k+1}$. Then

$$ck^p - 1 < \lfloor ck^p \rfloor \leq n \leq \lfloor c(k+1)^p \rfloor - 1 \leq c(k+1)^p - 1.$$

Solving $k = k_n$ from the preceding inequality, we obtain

(5)     $$t_n = a_{k_n}, \qquad \text{where } k_n = \lceil (c^{-1}(n+1))^{1/p} \rceil - 1 \quad \text{and}$$

$$\lceil u \rceil = \min\{i \in \mathbb{Z} : i \geq u\}.$$

With the above formula, it is easy to check the condition $n+1 = t_{n+1}$ in step 1 of Algorithm 1. In the special case with $c = 1$ and $p = 2$, $n+1 = t_{n+1}$ if and only if $(n+1)^{1/2} \in \mathbb{N}$.

Algorithm 1 is not yet directly applicable in practical situations since $\mu$ is unknown and $W_i$ needs to be centered. A natural centering sequence is the sample mean $\bar{X}_n = \sum_{i=1}^n X_i/n$. Based on $V_n$ in (3), we propose

(6)     $$V_n' = \sum_{i=1}^n (W_i')^2, \qquad \text{where } W_i' = X_{t_i} + X_{t_i+1} + \cdots + X_i - l_i \bar{X}_n,$$

where we recall $l_i = i - t_i + 1$. Observe that $(W_i')^2 - W_i^2 = (l_i \bar{X}_n)^2 - 2l_i W_i \bar{X}_n$. To recursively compute $V_n'$, we also need to introduce

$$U_n = \sum_{i=1}^n l_i W_i \quad \text{and} \quad q_n = \sum_{i=1}^n l_i^2.$$

Then

(7)     $$V_n' = V_n - 2U_n \bar{X}_n + q_n (\bar{X}_n)^2.$$

Algorithm 1 can be modified as follows:

ALGORITHM 2.   At stage $n$, we store $(n, k_n, a_{k_n}, v_n, q_n, U_n, V_n, W_n, \bar{X}_n)$. At stage $n+1$, we update the vector by:

1. $k_{n+1} = k_n + \mathbf{1}_{n+1=a_{1+k_n}}$, $t_{n+1} = a_{k_{n+1}}$,
2. $\bar{X}_{n+1} = (n\bar{X}_n + X_{n+1})/(n+1)$,
3. $q_{n+1} = q_n + (n+2 - t_{n+1})^2$,

4. $v_{n+1} = v_n + (n + 2 - t_{n+1})$,
5. $W_{n+1} = X_{n+1} + W_n \mathbf{1}_{n+1 \neq t_{n+1}}$,
6. $V_{n+1} = V_n + W_{n+1}^2$,
7. $U_{n+1} = U_n + (n + 2 - t_{n+1})W_{n+1}$,
8. $V'_{n+1} = V_{n+1} - 2U_{n+1}\bar{X}_{n+1} + q_{n+1}(\bar{X}_{n+1})^2$.

Output: $\hat{\sigma}_{n+1}^2 = V'_{n+1}/v_{n+1}$.

At stage $n$, based on $\hat{\sigma}_n^2 = V'_n/v_n$, we can construct the $(1 - \alpha)$ confidence interval for $\mu$ as $\bar{X}_n \pm \hat{\sigma}_n z_{1-\alpha/2}/\sqrt{n}$. Convergence rate of $\hat{\sigma}_n^2$ certainly depends on the sequence $(a_k)$, as well as the dependence structure of the underlying process. Section 3 concerns the convergence properties of $\hat{\sigma}_n^2$.

It is easily seen that the above recursive algorithms can be generalized to spectral density estimation. Let

$$f(\theta) = \frac{1}{2\pi} \sum_{k \in \mathbb{Z}} \gamma(k) e^{\sqrt{-1}k\theta} = \frac{1}{2\pi} \sum_{k \in \mathbb{Z}} \gamma(k) \cos(k\theta), \qquad \theta \in \mathbb{R},$$

be the spectral density function, where $\sqrt{-1}$ is the imaginary unit. Assume that $\mathbb{E}X_k = 0$. As in (3), we can introduce

$$V_n(\theta) = \sum_{i=1}^{n} |W_i(\theta)|^2, \qquad \text{where } W_i(\theta) = \sum_{j=t_i}^{i} X_j e^{\sqrt{-1}j\theta},$$

and recursively estimate $f(\theta)$ at a given $\theta \in \mathbb{R}$ by $\hat{f}_n(\theta) = V_n(\theta)/(2\pi v_n)$. The latter can be viewed as a version of Bartlett's spectral density estimate with varying block lengths. Using similar but lengthier arguments adopted in the Appendix, we can obtain similar convergence results for $\hat{f}_n(\theta)$. The details are omitted since our primary focus is the inference of sample means of stationary processes.

3. Convergence properties. For the recursive estimate $\hat{\sigma}_n^2$ proposed in Section 2, a natural question is to study its convergence properties. The latter problem is far from being trivial. Here we should implement the dependence measures proposed in Wu (2005) and obtain moment and almost sure convergence of $\hat{\sigma}_n^2$.

We first make some structural assumptions on the dependence. Assume hereafter that $(X_i)$ is a stationary causal process of the form

(8)                           $X_i = g(\ldots, \varepsilon_{i-1}, \varepsilon_i)$,

where $\varepsilon_i$ are i.i.d. innovations and $g$ is a measurable function such that $X_i$ is well defined. The framework (8) is very general and it allows many widely used linear and nonlinear processes. As in Wiener (1958) and Priestley (1988), (8) can be interpreted as a physical system with $\mathcal{F}_i = (\ldots, \varepsilon_{i-1}, \varepsilon_i)$ being the input, $g$ being a filter and $X_i$ being the output. Wiener (1958) dealt with the problem of representing stationary and ergodic processes as shifts of functions of independent random variables; see Rosenblatt (1959), Tong (1990) and Borkar (1993). Based on (8), Wu (2005) introduced the physical and predictive dependence measures which quantifies the degree of dependence of outputs on inputs. Specifically, let $\varepsilon'_0, \varepsilon_j, j \in \mathbb{Z}$, be i.i.d. random variables and $\mathcal{F}'_0 = (\ldots, \varepsilon_{-2}, \varepsilon_{-1}, \varepsilon'_0)$; let $g_i(\mathcal{F}_0) = \mathbb{E}[g(\mathcal{F}_i)|\mathcal{F}_0]$. For $p \geq 1$ define the physical dependence measure

$$(9) \qquad \delta_p(i) = \|X_i - X'_i\|_p \qquad \text{where } X'_i = g(\mathcal{F}'_0, \varepsilon_1, \ldots, \varepsilon_{i-1}, \varepsilon_i),$$

and the predictive dependence measure

$$(10) \qquad \omega_p(i) = \|g_i(\mathcal{F}_0) - g_i(\mathcal{F}'_0)\|_p.$$

The process $X'_i$ is a coupled version of $X_i$ with $\varepsilon_0$ replaced by $\varepsilon'_0$. So $\delta_p(i)$ quantifies the contribution of $\varepsilon_0$ to $X_i$ by measuring the distance between $X_i$ and $X'_i$. $\omega_p(i)$ measures the contribution of $\varepsilon_0$ in predicting future expected values. For details, see Wu (2005).

In comparison with the traditional strong mixing conditions, $\delta_p(i)$ and $\omega_p(i)$ appear more convenient to use in our context and they are directly related with the data-generating mechanisms. Wu (2005) showed that, if the process $(X_i)$ is stable, namely, $\Omega_2 := \sum_{i=0}^{\infty} \omega_2(i) < \infty$, then (13) below holds with $\sigma \leq \Omega_2$. See also Hannan (1979) and Volný (1993). Box, Jenkins and Reinsel (1994) considered the special case of linear processes and interpreted the stability condition as the cumulative impact of a single shock $\varepsilon_0$ on the whole process $(X_i)$ being finite. Main results in the sequel are all expressed in terms of $\delta_p(i)$ and $\omega_p(i)$.

3.1. *A representation of $\sigma$.* We shall first introduce a useful representation of $\sigma$. Write $S_i = \sum_{j=1}^{i} X_j$. Assume that $\mathbb{E}X_i = 0$ and

$$(11) \qquad \sum_{i=0}^{\infty} \|\mathcal{P}_0 X_i\|_2 < \infty \qquad \text{where } \mathcal{P}_i \cdot := \mathbb{E}(\cdot|\mathcal{F}_i) - \mathbb{E}(\cdot|\mathcal{F}_{i-1}).$$

Then

$$(12) \qquad D_k := \sum_{i=k}^{\infty} \mathcal{P}_k X_i \in \mathcal{L}^2$$

and $(D_k)_{k \in \mathbb{Z}}$ is a stationary martingale difference sequence with respect to the natural filter $\mathcal{F}_k$. Additionally, by Theorem 1 in Hannan (1979), we have the invariance principle

$$(13) \qquad \frac{1}{\sqrt{n}} \left\{ \sum_{i \leq nt} X_i, 0 \leq t \leq 1 \right\} \Rightarrow \{\sigma \mathbb{B}(t), 0 \leq t \leq 1\} \qquad \text{where } \sigma = \|D_k\|_2.$$

Here $\mathbb{B}$ is the standard Brownian motion. Let $M_n = \sum_{i=1}^n D_i$. If (11) holds with $\alpha > 2$ [cf. (14) below], then we have $\|S_n - M_n\|_\alpha = o(\sqrt{n})$ [see Theorem 1 in Wu (2007)]. The operator $\mathcal{P}_i$ in (14) is called the *projection operator* and it naturally generates martingale differences. The representation of $\sigma$ in (13) is useful in the analysis of our estimates.

3.2. *Moment convergence.* We first present a general result on moment convergence properties of $V_n/v_n$ under mild dependence conditions. Recall (11) for the definition of the projection operator $\mathcal{P}_i \cdot = \mathbb{E}(\cdot|\mathcal{F}_i) - \mathbb{E}(\cdot|\mathcal{F}_{i-1})$.

THEOREM 1. *Let $\mathbb{E}X_i = 0$ and $X_i \in \mathcal{L}^\alpha$, $\alpha > 2$. Assume*

$$\sum_{i=0}^\infty \|\mathcal{P}_0 X_i\|_\alpha < \infty. \tag{14}$$

*Further assume that, as $m \to \infty$, $a_{m+1} - a_m \to \infty$ and*

$$\frac{(a_{m+1} - a_m)^2}{\sum_{k=2}^m (a_k - a_{k-1})^2} \to 0. \tag{15}$$

*Then $\|V_n/v_n - \sigma^2\|_{\alpha/2} = o(1)$.*

Theorem 1 implies that, for consistency of $V_n/v_n$, $X_k$ does not need to have finite fourth moment. Instead, the moment condition $X_i \in \mathcal{L}^\alpha$ with $\alpha > 2$ suffices. We now discuss conditions (14) and (15) in the following remarks.

REMARK 1. By Jensen's inequality, we have $\|\mathcal{P}_0 X_i\|_\alpha \le \omega_\alpha(i) \le 2\|\mathcal{P}_0 X_i\|_\alpha$; see Theorem 1 in Wu (2005). Then (14) is equivalent to the stability condition $\sum_{j=0}^\infty \omega_\alpha(j) < \infty$ [Wu (2005)]. The latter condition can be interpreted as follows: the cumulative contribution of $\varepsilon_0$ in predicting future values $(X_i)_{i>0}$ is finite, thus suggesting short-range dependence. For long-range dependent processes (14) is violated and $\sigma^2$ does not always exist; see Example 5.2. So (14) is a very natural condition.

REMARK 2. Theorem 1 imposes mild conditions on the sequence $(a_k)_{k\ge 1}$. The theorem is applicable if $a_k = \lfloor ck^p \rfloor$, where $p > 1$ and $c > 0$ are constants. To account for dependence, it is certainly needed that $a_{m+1} - a_m \to \infty$. Condition (15) does not hold if $a_m$ diverges to infinity too fast. For example, (15) is violated if $a_k = 2^k$. In the latter case $V_n/v_n$ is *not* a consistent estimate of $\sigma^2$ if $X_i$ are i.i.d. standard normals. To see this, let $\xi_j$, $j \in \mathbb{Z}$, be independent and identically distributed as $\int_0^1 \mathbb{B}^2(t)\,dt$, where we recall that $\mathbb{B}$ is the standard Brownian motion. Elementary calculations show that $v_{2^m} \sim 2^{2m}/6$ and

$$2^{-2k} \sum_{i=2^k}^{2^{k+1}-1} (X_{2^k} + \cdots + X_i)^2 \Rightarrow \xi_0.$$

Since $X_i$ are i.i.d., $V_{2^m}/v_{2^m} \Rightarrow (3/2) \sum_{j=0}^\infty \xi_j/4^j$. In contrast, $\sigma = 1$.

Corollary 1 asserts the moment convergence of $\hat{\sigma}_n^2 = V'_n/v_n$ generated from Algorithm 2 which allows unknown $\mu$.

COROLLARY 1. *Let conditions (14) and (15) of Theorem 1 be satisfied. Then for $\hat{\sigma}_n^2 = V'_n/v_n$ generated from Algorithm 2, we also have $\|V'_n/v_n - \sigma^2\|_{\alpha/2} = o(1)$.*

3.3. *Convergence rates.* Theorem 1 asserts the moment convergence of $V_n/v_n$ under mild conditions (14) and (15). However, it does not provide information on the convergence rates. Under suitable decay rates of dependence measures, Theorem 2 provides a convergence rate of $V_n/v_n$ for algebraic sequences $(a_k)$.

THEOREM 2. *Let $a_k = \lfloor ck^p \rfloor$, $k \geq 1$, where $c > 0$ and $p > 1$ are constants.*

(i) *Assume that $X_i \in \mathcal{L}^\alpha$, $\mathbb{E}X_i = 0$, and for some $\alpha \in (2, 4]$,*

$$(16) \qquad\qquad \sum_{j=0}^\infty \delta_\alpha(j) < \infty.$$

*Then*

$$(17) \qquad\qquad \|V_n - \mathbb{E}V_n\|_{\alpha/2} = O(n^{3/2 - 3/(2p) + 2/\alpha}).$$

(ii) *Assume that $X_i \in \mathcal{L}^\alpha$, $\mathbb{E}X_i = 0$ and (16) holds for some $\alpha > 4$. Then*

$$(18) \qquad\qquad \lim_{n\to\infty} \frac{\|V_n - \mathbb{E}V_n\|}{n^{2-3/(2p)}} = \frac{\sigma^2 p^2 c^{3/(2p)}}{\sqrt{12p - 9}}.$$

(iii) *If $X_i \in \mathcal{L}^2$, $\mathbb{E}X_i = 0$, and for some $q \in (0, 1]$,*

$$(19) \qquad\qquad \sum_{j=0}^\infty j^q \omega(j) < \infty.$$

*Then $\mathbb{E}V_n - v_n\sigma^2 = O[n^{1+(1-q)(1-1/p)}]$. Consequently, under (16) and (19), $\|V_n - v_n\sigma^2\|_{\alpha/2} = O(n^\phi)$, where $\phi = \max(3/2 - 3/(2p) + 2/\alpha, 1 + (1-q)(1-1/p))$.*

Since $\omega(j) \leq \delta_2(j) \leq \delta_\alpha(j)$, a sufficient condition for (16) and (19) is $\sum_{j=1}^{\infty} j^q \delta_\alpha(j) < \infty$.

Theorem 2 gives guidance on how to choose $p$ based on the dependence and moment conditions of the process, which are characterized by parameters $q$ and $\alpha$, respectively. A good $p$ is the minimizer of $n^{3/2-3/(2p)+2/\alpha} + n^{1+(1-q)(1-1/p)}$. This $p$ also minimizes $\phi = \phi(p)$. Solving the equation

$$3/2 - 3/(2p) + 2/\alpha = 1 + (1-q)(1-1/p),$$

one obtains $p = (1/2+q)/(q-2+2/\alpha)$. To summarize, we have the following:

COROLLARY 2. Let $p = (1/2+q)/(q-1/2+2/\alpha)$. Under conditions of Theorem 2, we have $\|V_n/v_n - \sigma^2\|_{\alpha/2} = O(n^{2/\alpha-1/2-1/(2p)})$. In particular, if $\alpha = 4$ and $q = 1$, then $p = 3/2$ and $\|V_n/v_n - \sigma^2\|_2 = O(n^{-1/3})$.

REMARK 3. Since $a_{k+1} - a_k \sim cpk^{p-1}$ and $m \sim (n/c)^{1/p}$, elementary calculations show that

$$(20) \quad v_{a_m} \sim \sum_{i=1}^{m} (a_{i+1} - a_i)(a_{i+1} - a_i + 1)/2 \sim m^{2p-1}c^2p^2/(4p-2) \sim v_n.$$

By Remark 4, $\|V_n - V_n'\|_{\alpha/2}/v_n = O(n^{-1/p})$. Hence, Corollary 2 also applies to $\hat{\sigma}_n^2 = V_n'/v_n$ since $-1/p < 2/\alpha - 1/2 - 1/(2p)$.

Since $2 < \alpha \leq 4$, $p$ increases as $q$ decreases. The latter observation can be explained as follows: if (19) only holds for small $q$, then it indicates strong dependence and one needs to choose large block sizes to suppress the dependence.

We now compare Corollary 2 with classical results of the estimation of TAVC by using the batched means. Carlstein (1986) obtained the bound $O(n^{-1/3})$ for the special AR(1) model with i.i.d. normal innovations. Under appropriate strong mixing conditions, one can obtain the optimal mean squares error (MSE) bound $O(n^{-2/3})$ if the batch size is of order $n^{1/3}$; see Künsch (1989) and Lahiri (2003), among others. By Corollary 2, one can obtain the same bound: $\|V_n/v_n - \sigma^2\|_2^2 = O(n^{-2/3})$, and the gap $a_{m+1} - a_m = \lfloor c(m+1)^{3/2} \rfloor - \lfloor cm^{3/2} \rfloor \sim (3c/2)m^{1/2} \sim (c^{3/2}3/2)n^{1/3}$. For more discussions, see Section 4.

Our results have the attractive feature that they do not require strong mixing conditions which may be difficult to be verified in practice. Also, we impose a very mild moment condition that $X_i \in \mathcal{L}^\alpha$ with $2 < \alpha \leq 4$.

In view of the recursive nature of our estimate, it is natural to consider its almost sure convergence behavior. In the context of mean estimation based on MCMC simulations, Glynn and Whitt (1992) argued that, for asymptotic validity of sequential confidence intervals, one needs to have a strongly consistent estimate of $\sigma$ while the weaker version of mere convergence in probability is not enough.

COROLLARY 3. *Under conditions in Corollary* 2, *we have*

$$(21) \qquad \left\| \max_{n \leq N} |V_n - \mathbb{E}V_n| \right\|_{\alpha/2} = O(N^\tau \log N),$$

$$\text{where } \tau = 3/2 - 3/(2p) + 2/\alpha,$$

*and* $V_N - \mathbb{E}V_N = o[N^\tau(\log N)^2]$ *almost surely, and also*

$$(22) \qquad V_N/v_N - \sigma^2 = o(N^{2/\alpha - 1/2 - 1/(2p)}(\log N)^2) \qquad \text{almost surely.}$$

**4. Implementation issues.** Assume that (19) holds with $q = 1$ and (16) holds with $\alpha > 4$. Let the sequence $a_k = \lfloor ck^p \rfloor$, $k \geq 1$. To implement Algorithm 2, it is necessary to choose $c$ and $p$. Corollary 2 suggests the optimal $p = 3/2$. Here we shall suggest a data driven estimate of $c$ by using the procedure in Bühlmann and Künsch (1999).

Since (19) holds with $q = 1$, $\sum_{i=1}^\infty i|\gamma(i)| < \infty$. So as $l \to \infty$,

$$\mathbb{E}(S_l^2) - l\sigma^2 = -2\sum_{k=1}^\infty \min(l, k)\gamma(k) = \theta + o(1) \qquad \text{where } \theta = -2\sum_{k=1}^\infty k\gamma(k).$$

So $\mathbb{E}V_n - v_n\sigma^2 = n\theta + o(n)$. By (20), $v_n \sim 9m^2c^2/16$. Since $m \sim (n/c)^{2/3}$, by Theorem 2(ii),

$$\|V_n/v_n - \sigma^2\|_2^2 = \frac{\|V_n - \mathbb{E}V_n\|_2^2 + |\mathbb{E}V_n - v_n\sigma^2|^2}{v_n^2}$$

$$\sim \frac{16\sigma^4}{9m} + \frac{256\theta^2 n^2}{81c^4 m^4} \sim \left(\sigma^4\frac{16c^{2/3}}{9} + \theta^2\frac{256}{81c^{4/3}}\right)n^{-2/3}.$$

The MSE-optimal $c$ minimizes $\|V_n/v_n - \sigma^2\|_2^2$. Hence,

$$(23) \qquad \|V_n/v_n - \sigma^2\|_2^2 \sim \frac{2^{14/3}}{3^{5/3}}\theta^{2/3}\sigma^{8/3}n^{-2/3} \quad \text{and} \quad c = \frac{4\sqrt{2}|\theta|}{3\sigma^2}.$$

We now consider the batched mean estimate $\sigma_n^2(l_n)$ given in (2) with $\bar{X}_n$ therein replaced by 0. Assume $l_n/n \to 0$ and $l_n \to \infty$. Under suitable strong mixing conditions, we have $\|\sigma_n^2(l_n) - \mathbb{E}\sigma_n^2(l_n)\|_2^2 \sim 4\sigma^4 l_n/(3n)$ and $\mathbb{E}\sigma_n^2(l_n) - \sigma^2 \sim (\theta + o(1))/l_n$ [see, e.g., Song and Schmeiser (1995) or Politis, Romano and Wolf (1999)]. So the asymptotic MSE-optimal $l_n$ satisfies

$$(24) \qquad \|\sigma_n^2(l_n) - \sigma^2\|_2^2 \sim \frac{2^{2/3}3^{1/3}\theta^{2/3}\sigma^{8/3}}{n^{2/3}}$$

$$\text{with } l_n = \lfloor \lambda_* n^{1/3} \rfloor \text{ and } \lambda_*^3 = \frac{3\theta^2}{2\sigma^4}.$$

Bühlmann and Künsch ([1999](#)) proposed a data-driven method for finding the block length $l_n$. Sherman ([1998](#)) considered a similar problem. For the purpose of estimating $c$ in ([23](#)), we shall present Bühlmann and Künsch's ([1999](#)) algorithm here.

ALGORITHM 3. Let the Tukey–Hanning window $w_{TH}(x) = (1 + \cos(\pi x)) \times \mathbf{1}_{|x| \le 1}/2$ and the split-cosine window $w_{SC}(x) = (1 + \cos(5(x - 0.8)\pi))/2$ if $0.8 \le |x| \le 1$; $w_{SC}(x) = 1$ if $0.8 > |x|$ and $w_{SC}(x) = 0$ if $|x| > 1$.

1. Calculate $\hat{\gamma}(k) = n^{-1} \sum_{i=1}^{n-|k|} (X_i - \bar{X}_n)(X_{i+|k|} - \bar{X}_n)$, $k = 1 - n, \ldots, n - 1$.
2. Let $b_0 = n^{-1}$. For $m = 1, 2, 3, 4$, let

$$b_m = n^{-1/3} \left( \frac{\sum_{k=1-n}^{n-1} \hat{\gamma}(k)^2}{6 \sum_{k=1-n}^{n-1} w_{SC}(kb_{m-1}n^{4/21})k^2 \hat{\gamma}(k)^2} \right)^{1/3}.$$

3. Let $\hat{l}_n$ be the closest integer of $\hat{b}^{-1}$, where

$$\hat{b} = n^{-1/3} \left( \frac{2(\sum_{k=1-n}^{n-1} \hat{w}_{TH}(kb_4 n^{4/21}) \gamma(k))^2}{3(\sum_{k=1-n}^{n-1} w_{SC}(kb_4 n^{4/21})|k| \hat{\gamma}(k))^2} \right)^{1/3}.$$

By Theorem 4.1 in Bühlmann and Künsch's ([1999](#)), under suitable conditions, one has asymptotically that $n\hat{b}^3 \sim 2\sigma^4/(3\theta^2)$. Relation ([23](#)) hence suggests a data driven choice $\hat{c} = (4\hat{\lambda}_*/3)^{3/2}$, where $\hat{\lambda}_* = \hat{l}_n/n^{1/3}$ and $\hat{l}_n$ is from Algorithm [3](#). By ([23](#)) and ([24](#)), with $c = (4\lambda_*/3)^{3/2}$, we have $\|V_n/v_n - \sigma^2\|_2/\|\sigma_n^2(l_n) - \sigma^2\|_2 \sim 4/3$, which suggests that the recursive estimate $V_n/v_n$ has a reasonably good performance compared with the batched mean estimate $\sigma_n^2(l_n)$. In practice, we can conduct a pilot study and estimate $c$ by using Algorithm [3](#) with a relatively small $n$. Then we can use this $c$ for our recursive algorithm.

The computational and memory advantage of our recursive algorithm is more prominent if one runs multiple copies of the chain. In such applications we may obtain an estimate of $\sigma^2$ for each individual chain, and then use median or mean of those estimates to obtain an improved estimate. Also, we can check the variations of those TAVC estimates for convergence diagnostics. The computational cost for the traditional nonrecursive algorithms may be very expensive if the number of copies is large. Chauveau and Diebolt ([2003](#)) also considered estimate of $\sigma^2$ based on multiple chains. However, their estimate is not consistent if the number of copies is bounded.

5. **Applications.** Here we shall apply Theorems [1](#) and [2](#) to Markov chains which are in the form of iterated random functions and to functionals of linear processes. The former is useful in MCMC simulations.

5.1. *Markov chains.* Let $\varepsilon_i$, $i \in \mathbb{Z}$, be i.i.d. random variables. Consider the Markov chain $(Y_n)$ recursively defined by

$$(25) \qquad Y_n = g(Y_{n-1}, \varepsilon_n),$$

where $g$ is a measurable function. A variety of nonlinear time series models are of the form (25). Diaconis and Freedman (1999) showed that the Markov chain (25) admits a unique stationary distribution provided that

$$(26) \qquad \mathbb{E} \log L_{\varepsilon_0} < 0 \qquad \text{where } L_{\varepsilon_0} = \sup_{y \neq y'} \frac{|g(y, \varepsilon_0) - g(y', \varepsilon_0)|}{|y - y'|},$$

and

$$(27) \qquad \mathbb{E}[L_{\varepsilon_0}^\iota + |g(y_0, \varepsilon_0)|^\iota] < \infty \qquad \text{for some } y_0 \text{ and } \iota > 0.$$

Under (26) and (27), by iterating (25), $Y_n$ adopts the representation (8). Interestingly, the same set of conditions [namely, (26) and (27)] also implies that $\delta_\chi(j) = O(r^j)$ for some $r \in (0, 1)$ and $\chi > 0$; see Wu and Shao (2004).

We now apply Theorem 2 to the process $X_i = h(Y_i)$. In MCMC experiments, $\mu = \mathbb{E} X_i$ is estimated by $\bar{X}_n$ and the length of the confidence interval $\bar{X}_n \pm z_{1-\alpha/2} \hat{\sigma}_n / \sqrt{n}$ can be used for convergence diagnostics [Jones et al. (2006)]. We shall impose regularity conditions on $h$ such that (16) and (19) are satisfied. Assume $X_i \in \mathcal{L}^{\alpha_0}$ for some $\alpha_0 > 2$. For $t > 0$ let

$$\Delta(t) = \sup\{\|[h(Y) - h(Y')] \times \mathbf{1}_{|Y-Y'| \leq t}\|_\alpha : Y \text{ and } Y' \text{ are identically distributed}\}.$$

Following the argument of Theorem 3 in Wu and Shao (2004), under

$$(28) \qquad \int_0^1 \frac{\Delta(t) |\log t|}{t} \, dt < \infty,$$

we have $\sum_{i=1}^\infty i \delta_\alpha(i) < \infty$ and, hence, (16) and (19) hold. The details of the derivation are omitted. We now give examples that (28) holds. If $h$ is Lipschitz continuous, then $\Delta(t) = O(t)$ and (28) follows. Let $h$ be an indicator function $h(y) = \mathbf{1}_{y \leq y_0}$, where $y_0$ is fixed. Then (28) also holds if $\mathbb{P}(|Y_i - y_0| \leq t) = O(t^\rho)$ for some $\rho > 0$. In particular, if $Y_i$ has a density, then $\rho = 1$.

An attractive feature of our setting is that we do not need the assumption of irreducibility and positive Harris recurrence. The latter assumptions are not valid for many Markov chains. For example, Markov chains associated with fractal images [Diaconis and Freedman (1999)] are not generally positive Harris recurrent. As a concrete example, consider (25) with $Y_n = (Y_{n-1} + 2\varepsilon_n)/3$, where $\varepsilon_n$ are i.i.d. with distribution $\mathbb{P}(\varepsilon_n = 0) = \mathbb{P}(\varepsilon_n = 1) = 1/2$. Then the chain is not positive Harris recurrent. On the other hand, (26) and (27) trivially hold and $(Y_n)$ adopts an invariant distribution. Additionally, its support is the Cantor set and $\mathbb{P}(|Y_i - y_0| \leq t) = O(t^\rho)$, where $\rho = (\log 2)/(\log 3)$ is the Hausdorff dimension.

5.2. *Linear processes.* Let $\varepsilon_i$, $i \in \mathbb{Z}$, be i.i.d. random variables with mean 0 and finite $\alpha$th moment ($\alpha > 2$) and $(a_i)$ be a sequence of real coefficients; let $X_n = K(e_n)$, where $K$ is a measure function for which $X_n \in \mathcal{L}^\alpha$ and $e_n = \sum_{i=0}^{\infty} a_i \varepsilon_{n-i}$ is a linear process. A special case is that $K(x) = |x|$. Since $K$ may be nonlinear, the treatment of $\sum_{i=1}^{n} K(e_i)$ appears more difficult than that of $\sum_{i=1}^{n} e_i$ since the latter preserves the linearity structure.

We now apply Theorem 1 to the process $(X_i)$. Recall that $\varepsilon_0'$ is independent of $\varepsilon_i$, $i \in \mathbb{Z}$. Let $e_n' = e_n - a_n \varepsilon_0 + a_n \varepsilon_0'$. If $K$ is Lipschitz continuous, then $|K(e_n) - K(e_n')| = O(|a_n|)|\varepsilon_0 - \varepsilon_0'|$. Hence, the physical dependence measure $\delta_\alpha(n) = O(|a_n|)$ and, consequently, $\|\mathcal{P}_0 X_i\|_\alpha = O(|a_i|)$ since $\omega_\alpha(n) \leq \delta_\alpha(n)$. In this case (14) is reduced to $\sum_{i=0}^{\infty} |a_i| < \infty$, which is a natural condition for the short-range dependence. If the latter condition is violated, for example, if $a_i = i^{-\beta}$, $1/2 < \beta < 1$, then the $(X_i)$ is a long-memory process and normalizing sequence for $\sum_{i=1}^{n} X_i$ is $n^{3/2-\beta}$, which is different from $\sqrt{n}$. Correspondingly, $\sigma^2 = \infty$.

## APPENDIX

**A.1. Proof of Theorem 1.** For $n \in \mathbb{N}$ choose $m = m_n \in \mathbb{N}$ such that $a_m \leq n < a_{m+1}$. Then

$$(29) \quad \begin{aligned} v_n &= \sum_{j=1}^{n} (j - t_j + 1) = \sum_{i=2}^{m} \sum_{j=a_{i-1}}^{a_i - 1} (j - t_j + 1) + \sum_{j=a_m}^{n} (j - t_j + 1) \\ &= \sum_{i=2}^{m} \frac{(a_i - a_{i-1})(a_i - a_{i-1} + 1)}{2} + \frac{(n - a_m)(n - a_m + 1)}{2}. \end{aligned}$$

Simple calculations show that (15) implies

$$(30) \quad 1 \leq \liminf_{m \to \infty} \frac{v_n}{v_{a_m}} \leq \limsup_{m \to \infty} \frac{v_{a_{m+1}}}{v_{a_m}} = 1.$$

So the limits in the above expression are all 1. Also observe that for any fixed $k_0 \in \mathbb{N}$, since $a_{m+1} - a_m$ is increasing to $\infty$, we have

$$(31) \quad \lim_{m \to \infty} \frac{\#\{i \leq n : i - t_i + 1 \leq k_0\}}{v_n} \leq \lim_{m \to \infty} \frac{m k_0}{v_n} = 0.$$

We now apply the martingale approximation in Wu (2007). Clearly (14) implies that $D_k := \sum_{i=k}^{\infty} \mathcal{P}_k X_i \in \mathcal{L}^\alpha$. Let $M_n = \sum_{i=1}^{n} D_i$. By Theorem 1 in Wu (2007), condition (14) also implies that

$$(32) \quad \|S_n\|_\alpha = O(\sqrt{n}), \qquad \|M_n\|_\alpha = O(\sqrt{n}) \quad \text{and} \quad \|S_n - M_n\|_\alpha = o(\sqrt{n}).$$

Hence, as $n \to \infty$,

$$(33) \quad \rho_n := n^{-1}\|S_n^2 - M_n^2\|_{\alpha/2} \leq n^{-1}\|S_n - M_n\|_\alpha \|S_n + M_n\|_\alpha \to 0.$$

As $V_n$ in (3), we introduce

$$Q_n = \sum_{i=1}^{n} R_i^2 \qquad \text{where } R_i = D_{t_i} + D_{t_i+1} + \cdots + D_i.$$

Our plan is to first approximate $V_n$ by $Q_n$ such that $\|Q_n - V_n\|_{\alpha/2} = o(v_n)$ and then show that $\|Q_n/v_n - \sigma^2\|_{\alpha/2} = o(1)$. Clearly the theorem follows from these two assertions. For the former, let $k_0 \in \mathbb{N}$. By (33) and (31),

$$
\begin{aligned}
(34) \qquad \limsup_{n\to\infty} \frac{\|V_n - Q_n\|_{\alpha/2}}{v_n} &\leq \limsup_{n\to\infty} v_n^{-1} \sum_{i=1}^{n} \|R_i^2 - W_i^2\|_{\alpha/2} \\
&\leq \limsup_{n\to\infty} v_n^{-1} \sum_{i=1}^{n} (i - t_i + 1)\rho_{i-t_i+1} \\
&\leq \limsup_{n\to\infty} v_n^{-1} \sum_{1 \leq i \leq n \,:\, i-t_i+1 > k_0} (i - t_i + 1)\rho_{i-t_i+1} \\
&\leq \sup_{k \geq k_0} \rho_k \to 0 \qquad \text{as } k_0 \to \infty.
\end{aligned}
$$

It remains to prove $\|Q_n/v_n - \sigma^2\|_{\alpha/2} = o(1)$. Note that $t_i = a_k$ if $a_k \leq i \leq a_{k+1} - 1$. Let

$$Y_k = \sum_{i=a_k}^{a_{k+1}-1} (D_{t_i} + D_{t_i+1} + \cdots + D_i)^2 = \sum_{i=a_k}^{a_{k+1}-1} (D_{a_k} + D_{a_k+1} + \cdots + D_i)^2$$

and

$$\tilde{Y}_k = \sum_{i=a_k}^{a_{k+1}-1} (D_{a_k}^2 + D_{a_k+1}^2 + \cdots + D_i^2).$$

By Burkholder's inequality, there exists a constant $c = c_\alpha$ such that

$$
\begin{aligned}
\|Y_k\|_{\alpha/2} &\leq \sum_{i=a_k}^{a_{k+1}-1} \|(D_{a_k} + D_{a_k+1} + \cdots + D_i)^2\|_{\alpha/2} \\
&= \sum_{i=a_k}^{a_{k+1}-1} \|D_{a_k} + D_{a_k+1} + \cdots + D_i\|_{\alpha}^2 \\
&\leq \sum_{i=a_k}^{a_{k+1}-1} c_\alpha (i - a_k + 1)\|D_1\|_{\alpha}^2.
\end{aligned}
$$

On the other hand,

$$\|\tilde{Y}_k\|_{\alpha/2} \leq \sum_{i=a_k}^{a_{k+1}-1} (i - a_k + 1)\|D_1\|_{\alpha}^2.$$

In the rest of the proof, $c_\alpha$ denotes a constant which only depends on $\alpha$ and its value may change from line to line. Since $1 < \alpha/2 \leq 2$ and $Y_k - \mathbb{E}(Y_k | \mathcal{F}_{a_k})$, $k = 1, 2, \ldots$, is a martingale difference sequence, we have by Burkholder's and Jensen's inequalities that

$$(35) \quad \begin{aligned} \left\| \sum_{k=1}^m [Y_k - \mathbb{E}(Y_k | \mathcal{F}_{a_k})] \right\|_{\alpha/2}^{\alpha/2} &\leq c_\alpha \sum_{k=1}^m \| Y_k - \mathbb{E}(Y_k | \mathcal{F}_{a_k}) \|_{\alpha/2}^{\alpha/2} \\ &\leq c_\alpha \sum_{k=1}^m \| Y_k \|_{\alpha/2}^{\alpha/2}. \end{aligned}$$

Similarly,

$$(36) \quad \left\| \sum_{k=1}^m [\tilde{Y}_k - \mathbb{E}(\tilde{Y}_k | \mathcal{F}_{a_k})] \right\|_{\alpha/2}^{\alpha/2} \leq c_\alpha \sum_{k=1}^m \| \tilde{Y}_k \|_{\alpha/2}^{\alpha/2}.$$

Note that $D_i$ are also martingale differences. Simple calculations show that $\mathbb{E}(\tilde{Y}_k | \mathcal{F}_{a_k}) = \mathbb{E}(Y_k | \mathcal{F}_{a_k})$. By (35) and (36),

$$\left\| \sum_{k=1}^m (Y_k - \tilde{Y}_k) \right\|_{\alpha/2}^{\alpha/2}$$

$$\leq c_\alpha \sum_{k=1}^m (\| Y_k \|_{\alpha/2}^{\alpha/2} + \| \tilde{Y}_k \|_{\alpha/2}^{\alpha/2})$$

$$\leq c_\alpha \| D_1 \|_\alpha^\alpha \sum_{k=1}^m \left[ \sum_{i=a_k}^{a_{k+1}-1} (i - a_k + 1) \right]^{\alpha/2}$$

$$\leq c_\alpha \| D_1 \|_\alpha^\alpha \max_{h \leq m} \left[ \sum_{i=a_h}^{a_{h+1}-1} (i - a_h + 1) \right]^{\alpha/2 - 1} \sum_{k=1}^m \left[ \sum_{i=a_k}^{a_{k+1}-1} (i - a_k + 1) \right].$$

By (15) and (30), since $a_{h+1} - a_h \to \infty$,

$$(37) \quad \frac{\| \sum_{k=1}^m (Y_k - \tilde{Y}_k) \|_{\alpha/2}^{\alpha/2}}{v_n^{\alpha/2}} \leq c_\alpha \| D_1 \|_\alpha^\alpha \left[ \frac{\max_{h \leq m} (a_{h+1} - a_h)}{v_n} \right]^{\alpha/2 - 1} \to 0.$$

By the ergodic theorem, since $D_k^2 \in \mathcal{L}^{\alpha/2}$, we have $\| D_1^2 + \cdots + D_l^2 - l\sigma^2 \|_{\alpha/2} = o(l)$. Therefore, $\| \tilde{Y}_k - \mathbb{E}\tilde{Y}_k \|_{\alpha/2} = o((a_{k+1} - a_k)^2)$ and, by (35) and (36),

$$\lim_{n \to \infty} \frac{\| \sum_{k=1}^m (\tilde{Y}_k - \mathbb{E}\tilde{Y}_k) \|_{\alpha/2}}{v_n} = \lim_{n \to \infty} \frac{\sum_{k=1}^m o((a_{k+1} - a_k)^2)}{v_n} = 0,$$

which, in view of (37), implies that $\| \sum_{k=1}^m Y_k - v_{a_m} \sigma^2 \|_{\alpha/2} = o(v_{a_m})$.

Finally, we shall compare $Q_n$ and $Q_{a_{m+1}-1} = \sum_{k=1}^{m} Y_k$. To this end, again by (35) and (36), recall $a_m \le n < a_{m+1}$,

$$\|Q_n - Q_{a_{m+1}-1}\|_{\alpha/2} = \left\| \sum_{i=n+1}^{a_{m+1}-1} R_i^2 \right\|_{\alpha/2} \le \sum_{i=n+1}^{a_{m+1}-1} \|R_i\|_\alpha^2$$

$$= \sum_{i=n+1}^{a_{m+1}-1} O(i - t_i + 1) \le (a_{m+1} - a_m)^2 = o(v_n),$$

which by (34) completes the proof.

**A.2. Proof of Corollary 1.** Observe that $V_n'$ remains unchanged if $X_i$ is replaced by $X_i - \mu$. So we can assume without loss of generality that $\mu = 0$. By (7) and Theorem 1, it suffices to verify that (i) $\|U_n \bar{X}_n\|_{\alpha/2} = o(v_n)$ and (ii) $\|q_n(\bar{X}_n)^2\|_{\alpha/2} = o(v_n)$. For (ii), by (32), $\|\bar{X}_n\|_\alpha = O(n^{-1/2})$. Choose $m \in \mathbb{N}$ such that $a_m \le n < a_{m+1}$. By (15),

$$(a_{m+1} - a_m)^2 = o(1) \left[ \sum_{k=2}^{m} (a_k - a_{k-1}) \right]^2 = o(a_m^2).$$

Since $a_m \to \infty$ and $a_m$ is increasing,

(38)  $$\max_{l \le m} (a_{l+1} - a_l) = o(a_m) = o(n).$$

Hence, $q_n \le v_n \max_{l \le m}(a_{l+1} - a_l) = v_n o(n)$ and (ii) follows. To show (i), we claim that

(39)  $$\|U_n\|_\alpha = O(1) \left[ \sum_{l=1}^{m} (a_{l+1} - a_l)^5 \right]^{1/2}.$$

With the above relation, noting that $\sum_{l=1}^{m}(a_{l+1} - a_l)^4 \le [\sum_{l=1}^{m}(a_{l+1} - a_l)^2]^2$, we have by (29) and (38) that $\|U_n \bar{X}_n\|_{\alpha/2} \le \|U_n\|_\alpha \|\bar{X}_n\|_\alpha = o(v_n)$.

In the sequel we shall prove (39). To this end, recall $l_i = i - t_i + 1$ and let

$$h_j = h_{j,n} = \sum_{i=1}^{n} l_i \mathbf{1}_{t_i \le j \le i}, \qquad j = 1, \dots, n.$$

Then

$$U_n = \sum_{i=1}^{n} l_i \sum_{j=t_i}^{i} X_j = \sum_{j=1}^{n} X_j h_j.$$

Since $X_j = \sum_{k=0}^{\infty} \mathcal{P}_{j-k} X_j$, and $\mathcal{P}_{j-k} X_j$, $j \in \mathbb{Z}$, forming martingale differences, we have by Burkholder's and Minkowski's inequalities that

$$\|U_n\|_\alpha \le \sum_{k=0}^{\infty} \left\| \sum_{j=1}^{n} \mathcal{P}_{j-k} X_j h_j \right\|_\alpha$$

$$\leq \sum_{k=0}^{\infty} c_\alpha \left[ \sum_{j=1}^{n} \|\mathcal{P}_{j-k} X_j h_j\|_\alpha^2 \right]^{1/2}$$

$$= \left( \sum_{j=1}^{n} h_j^2 \right)^{1/2} c_\alpha \sum_{k=0}^{\infty} \|\mathcal{P}_0 X_k\|_\alpha.$$

By (14) and the definition of $h_j$, (39) follows from

$$\sum_{j=1}^{n} h_j^2 \leq \sum_{k=1}^{m} \sum_{j=a_k}^{a_{k+1}-1} h_j^2 \leq \sum_{k=1}^{m} \sum_{j=a_k}^{a_{k+1}-1} (a_{k+1} - a_k)^4 = \sum_{k=1}^{m} (a_{k+1} - a_k)^5.$$

REMARK 4. If $a_k = \lfloor ck^p \rfloor$, $k \geq 1$, where $c > 0$ and $p > 1$, then $m \sim (n/c)^{1/p}$ and, by (39), $\|U_n\|_\alpha = O[m^{(5p-4)/2}] = O(n^{5/2-2/p})$. Also note that $q_n \asymp n^{3-2/p}$. Hence, $\|V_n - V_n'\|_{\alpha/2} = O(q_n/n) + O(n^{5/2-2/p})/n^{1/2} = O(n^{2-2/p})$.

**A.3. Proof of Theorem 2.** (i) Recall (3) for $W_i = X_{t_i} + X_{t_i+1} + \cdots + X_i$ and (9) for the definition of the coupled process $(X_n')$. Let $W_i^* = X_{t_i}' + X_{t_i+1}' + \cdots + X_i'$. For notational simplicity write $\delta_j$ for $\delta_\alpha(j)$. Since $\varepsilon_0'$ is independent of $\varepsilon_i, i \in \mathbb{Z}$, we have $\mathbb{E}(X_i|\mathcal{F}_{-1}) = \mathbb{E}(X_i^*|\mathcal{F}_{-1}) = \mathbb{E}(X_i^*|\mathcal{F}_0)$. By Jensen's inequality, $\|\mathcal{P}_0 X_i\|_\alpha \leq \|X_i - X_i^*\|_\alpha = \delta_i$ and (16) implies that $\Theta_\alpha = \sum_{i=0}^{\infty} \|\mathcal{P}_0 X_i\|_\alpha < \infty$. By Theorem 1 in Wu (2007), $\|W_i\|_\alpha \leq c_\alpha \Theta_\alpha (i - t_i + 1)^{1/2}$, where $c_\alpha$ is a constant. Since

$$\mathbb{E}[W_i^2|\mathcal{F}_{-1}] = \mathbb{E}[(W_i^*)^2|\mathcal{F}_{-1}] = \mathbb{E}[(W_i^*)^2|\mathcal{F}_0],$$

we have by Schwarz's and Jensen's inequalities that

$$\|\mathcal{P}_0 W_i^2\|_{\alpha/2} = \|\mathbb{E}[W_i^2|\mathcal{F}_0] - \mathbb{E}[W_i^2|\mathcal{F}_{-1}]\|_{\alpha/2}$$

$$= \|\mathbb{E}[W_i^2|\mathcal{F}_0] - \mathbb{E}[(W_i^*)^2|\mathcal{F}_0]\|_{\alpha/2}$$

$$\leq \|W_i^2 - (W_i^*)^2\|_{\alpha/2} \leq \|W_i + W_i^*\|_\alpha \|W_i - W_i^*\|_\alpha$$

$$\leq 2\|W_i\|_\alpha \sum_{j=t_i}^{i} \delta_j \leq 2c_\alpha \Theta_\alpha (i - t_i + 1)^{1/2} \sum_{j=t_i}^{i} \delta_j.$$

Similarly, for $k \geq 0$,

$$\|\mathcal{P}_{i-k} W_i^2\|_{\alpha/2} \leq 2\|W_i\|_\alpha \sum_{j=t_i}^{i} \delta_{k+t_i-j}$$

(40)

$$\leq 2c_\alpha \Theta_\alpha (i - t_i + 1)^{1/2} \sum_{j=t_i}^{i} \delta_{k+t_i-j}.$$

Since $\mathcal{P}_{i-k}W_i^2$, $i \in \mathbb{Z}$, form martingale differences, by Burkholder's inequality,

$$\left\|\sum_{i=1}^{n}\mathcal{P}_{i-k}W_i^2\right\|_{\alpha/2}^{\alpha/2} \leq c_\alpha \sum_{i=1}^{n}\|\mathcal{P}_{i-k}W_i^2\|_{\alpha/2}^{\alpha/2}$$

$$\leq c_\alpha \Theta_\alpha^{\alpha/2}\sum_{i=1}^{n}\left[(i-t_i+1)^{1/2}\sum_{j=t_i}^{i}\delta_{k+t_i-j}\right]^{\alpha/2}.$$

By the triangle inequality, since $W_i^2 = \sum_{k=0}^{\infty}\mathcal{P}_{i-k}W_i^2$, we have

$$(41) \qquad \|V_n - \mathbb{E}V_n\|_{\alpha/2} \leq \sum_{k=0}^{\infty}\left\|\sum_{i=1}^{n}\mathcal{P}_{i-k}W_i^2\right\|_{\alpha/2}.$$

If $a_m \leq i < a_{m+1}$, then $t_i = a_m$ and $i - t_i \leq a_{m+1} - 1 - a_m$. Let $b_m = \lfloor(1+c)p2^p m^{p-1}\rfloor$. Elementary calculations show that $a_{m+1} - 1 - a_m \leq b_m$ for all $m \in \mathbb{N}$. Hence,

$$\sum_{k=2b_m}^{\infty}\left\|\sum_{i=1}^{n}\mathcal{P}_{i-k}W_i^2\right\|_{\alpha/2}$$

$$(42) \qquad \leq \sum_{k=2b_m}^{\infty}\left\{\sum_{i=1}^{n}\left[(i-t_i+1)^{1/2}\sum_{j=0}^{b_m}\delta_{k-j}\right]^{\alpha/2}\right\}^{2/\alpha}O(1)$$

$$\leq \left[\sum_{i=1}^{n}(i-t_i+1)^{\alpha/4}\right]^{2/\alpha}\sum_{k=2b_m}^{\infty}\sum_{j=0}^{b_m}\delta_{k-j}O(1)$$

$$= [O(nb_m^{\alpha/4})]^{2/\alpha}o(b_m) = o(n^{2/\alpha}b_m^{3/2}).$$

On the other hand,

$$(43) \qquad \sum_{k=0}^{2b_m-1}\left\|\sum_{i=1}^{n}\mathcal{P}_{i-k}W_i^2\right\|_{\alpha/2} = O(b_m)\left[\sum_{i=1}^{n}(i-t_i+1)^{\alpha/4}\right]^{2/\alpha} = O(n^{2/\alpha}b_m^{3/2}).$$

Therefore, $\|V_n - \mathbb{E}V_n\|_{\alpha/2} = O(n^{2/\alpha}b_m^{3/2})$ and (i) follows since $b_m = O(n^{1-1/p})$.

(ii) Define $G_{h+1} = \sum_{i=a_h}^{a_{h+1}-1}W_i^2$. By Lemma 1 below, we have $\|G_{h+1} - \mathbb{E}(G_{h+1}|\mathcal{F}_{a_h})\|^2/(a_{h+1}-a_h)^4 \to \sigma^4/3$ as $h \to \infty$. Since $G_{h+1} - \mathbb{E}(G_{h+1}|\mathcal{F}_{a_h})$, $h = 1, 2, \ldots$, are martingale differences with respect to the filter $\mathcal{F}_{a_{h+1}}$, we have

$$\left\|\sum_{h=1}^{m}[G_{h+1} - \mathbb{E}(G_{h+1}|\mathcal{F}_{a_h})]\right\|^2 = \sum_{h=1}^{m}\|G_{h+1} - \mathbb{E}(G_{h+1}|\mathcal{F}_{a_h})\|^2$$

$$\sim \sum_{h=1}^{m}(a_{h+1}-a_h)^4\frac{\sigma^4}{3} \sim n^{4-3/p}\frac{\sigma^4 p^4 c^{3/p}}{12p-9}$$

by noting that $a_{h+1} - a_h \sim cph^{p-1}$. Similarly, by Lemma 1,

$$\left\| \sum_{h=1}^{m} [\mathbb{E}(G_{h+1}|\mathcal{F}_{a_h}) - \mathbb{E}(G_{h+1}|\mathcal{F}_{a_{h-1}})] \right\|^2$$

$$= \sum_{h=1}^{m} \| [\mathbb{E}(G_{h+1}|\mathcal{F}_{a_h}) - \mathbb{E}(G_{h+1}|\mathcal{F}_{a_{h-1}})] \|^2$$

$$\leq \sum_{h=1}^{m} \| \mathbb{E}(G_{h+1}|\mathcal{F}_{a_h}) - \mathbb{E}(G_{h+1}) \|^2$$

$$= \sum_{h=1}^{m} o((a_{h+1} - a_h)^4) = o(n^{4-3/p}).$$

We now deal with $\Xi_m := \sum_{h=1}^{m} [\mathbb{E}(G_{h+1}|\mathcal{F}_{a_{h-1}}) - \mathbb{E}(G_{h+1})]$. For $a_h \leq i \leq a_{h+1} - 1$, since $\mathbb{E}(W_i^2|\mathcal{F}_{a_{h-1}}) - \mathbb{E}(W_i^2) = \sum_{k=0}^{\infty} \mathcal{P}_{i-k} \mathbb{E}(W_i^2|\mathcal{F}_{a_{h-1}})$, we have

$$\|\Xi_m\| \leq \sum_{k=0}^{\infty} \left\| \sum_{h=1}^{m} \sum_{i=a_h}^{a_{h+1}-1} \mathcal{P}_{i-k} \mathbb{E}(W_i^2|\mathcal{F}_{a_{h-1}}) \right\|$$

$$= \sum_{k=0}^{\infty} \left[ \sum_{h=1}^{m} \sum_{i=a_h}^{a_{h+1}-1} \| \mathcal{P}_{i-k} \mathbb{E}(W_i^2|\mathcal{F}_{a_{h-1}}) \|^2 \right]^{1/2}.$$

Observe that $\mathcal{P}_{i-k} \mathbb{E}(W_i^2|\mathcal{F}_{a_{h-1}}) = 0$ if $i - k > a_{h-1}$, and $\mathcal{P}_{i-k} \mathbb{E}(W_i^2|\mathcal{F}_{a_{h-1}}) = \mathcal{P}_{i-k} W_i^2$ if $i - k \leq a_{h-1}$. By (40), as in the proof of (42), we have

$$\sum_{k=2b_m}^{\infty} \left[ \sum_{h=1}^{m} \sum_{i=a_h}^{a_{h+1}-1} \| \mathcal{P}_{i-k} \mathbb{E}(W_i^2|\mathcal{F}_{a_{h-1}}) \|^2 \right]^{1/2} = o(n^{1/2} b_m^{3/2}).$$

For $0 \leq k \leq 2b_m - 1$, since $\Theta_\alpha(l) = \sum_{i=l}^{\infty} \delta_\alpha(i) \to 0$ as $l \to \infty$,

$$\sum_{h=1}^{m} \sum_{i=a_h}^{a_{h+1}-1} \| \mathcal{P}_{i-k} \mathbb{E}(W_i^2|\mathcal{F}_{a_{h-1}}) \|^2$$

$$= O(1) \sum_{h=1}^{m} \sum_{i=a_h}^{a_{h+1}-1} (i - t_i + 1) \left[ \sum_{j=k+t_i-i}^{k} \delta_\alpha(j) \right]^2 \mathbf{1}_{i-k \leq a_{h-1}}$$

$$= O(1) \sum_{h=1}^{m} \sum_{i=a_h}^{a_{h+1}-1} (i - t_i + 1) \Theta_\alpha^2(a_h - a_{h-1})$$

$$= O(1) \sum_{h=1}^{m} (a_{h+1} - a_h)^2 \Theta_\alpha^2(a_h - a_{h-1})$$

$$= \sum_{h=1}^{m} o(h^{2p-2}) = o(m^{2p-1}).$$

Hence,

$$\sum_{k=0}^{2b_m-1} \left[ \sum_{h=1}^{m} \sum_{i=a_h}^{a_{h+1}-1} \|\mathcal{P}_{i-k}\mathbb{E}(W_i^2|\mathcal{F}_{a_{h-1}})\|^2 \right]^{1/2} = o(b_m m^{p-1/2})$$

and (ii) follows in view of

$$\left\| \sum_{i=n}^{a_{m+1}-1} (W_i^2 - \mathbb{E}(W_i^2)) \right\| \leq \sum_{i=n}^{a_{m+1}-1} \|W_i^2\| = O(b_m^2) = o(b_m m^{p-1/2})$$

since $\|W_i\|_4^2 = O(i - t_i + 1) = O(b_m)$, $a_{m+1} - 1 - a_m \leq b_m$ and $b_m = \lfloor (1 + c)p 2^p m^{p-1} \rfloor$.

(iii) Let $j > 0$. For $i \in \mathbb{Z}$, since $\mathcal{P}_i$ are orthogonal and $X_j = \sum_{i \in \mathbb{Z}} \mathcal{P}_i X_j$,

$$|\gamma(j)| = |\mathbb{E}(X_0 X_j)| = \left| \mathbb{E} \sum_{i \in \mathbb{Z}} (\mathcal{P}_i X_0)(\mathcal{P}_i X_j) \right|$$

$$\leq \sum_{i \in \mathbb{Z}} \|\mathcal{P}_i X_0\| \|\mathcal{P}_i X_j\| \leq \sum_{i \in \mathbb{Z}} \omega(-i)\omega(j-i).$$

Here we let $\omega(i) = 0$ if $i < 0$. By (19),

$$\sum_{j=0}^{\infty} j^q |\gamma(j)| < \infty.$$

Consequently, for $S_l = X_1 + \cdots + X_l$, since $0 < q \leq 1$,

$$|\mathbb{E}S_l^2 - l\sigma^2| \leq 2 \sum_{j=1}^{\infty} \min(j,l)|\gamma(j)| = O(l^{1-q}).$$

Therefore,

$$|\mathbb{E}V_n - t_n\sigma^2| \leq \sum_{i=1}^{n} |\mathbb{E}W_i - (i - t_i + 1)\sigma^2|$$

$$= \sum_{i=1}^{n} O[(i - t_i + 1)^{1-q}]$$

$$= O(nb_m^{1-q}) = O[n^{1+(1-q)(1-1/p)}].$$

LEMMA 1. *Assume that $X_i \in \mathcal{L}^\alpha$, $\mathbb{E}X_i = 0$ and (16) holds for some $\alpha > 4$. Let $S_i = \sum_{j=1}^i X_j$. Then we have* (i) $\|\sum_{i=1}^l (\mathbb{E}(S_i^2|\mathcal{F}_1) - \mathbb{E}(S_i^2))\| = o(l^2)$ *and* (ii)

$$\lim_{l \to \infty} \frac{\|\sum_{i=1}^l (S_i^2 - \mathbb{E}(S_i^2))\|^2}{l^4} = \frac{\sigma^4}{3}. \tag{44}$$

PROOF. As in (40), for $r \le 1$, $\|\mathcal{P}_r S_i^2\| \le C i^{1/2} \sum_{j=1}^i \delta_\alpha(j-r)$, where $C = 2c_\alpha \Theta_\alpha$. Since $\sum_{i=1}^l (\mathbb{E}(S_i^2|\mathcal{F}_1) - \mathbb{E}(S_i^2)) = \sum_{r=-\infty}^1 \sum_{i=1}^l \mathcal{P}_r S_i^2$, by orthogonality, (i) follows from

$$
\begin{aligned}
\left\| \sum_{i=1}^l (\mathbb{E}(S_i^2|\mathcal{F}_1) - \mathbb{E}(S_i^2)) \right\|^2 &= \sum_{r=-\infty}^1 \left\| \sum_{i=1}^l \mathcal{P}_r S_i^2 \right\|^2 \\
&\le \sum_{r=-\infty}^1 \left( \sum_{i=1}^l \|\mathcal{P}_r S_i^2\| \right)^2 \\
&\le \sum_{r=-\infty}^1 \left( C l^{3/2} \sum_{j=1}^l \delta_\alpha(j-r) \right)^2 \\
&\le \sum_{r=-\infty}^1 C^2 l^3 \Theta_\alpha \sum_{j=1}^l \delta_\alpha(j-r) \\
&= O(l^3) \sum_{j=1}^l \sum_{r=-\infty}^1 \delta_\alpha(j-r) = o(l^4).
\end{aligned}
$$

For (ii), let $A_l = \sum_{i=1}^l S_i^2/l^2$. By the invariance principle (13) and the continuous mapping theorem, we have $A_l \Rightarrow \sigma^2 \int_0^1 IB(t)^2 \, dt$. By Theorem 1 in Wu (2007), $\|S_i\|_\alpha = O(\sqrt{i})$. So

$$\|A_l\|_{\alpha/2} \le \sum_{i=1}^l \frac{\|S_i^2\|_{\alpha/2}}{l^2} \le \sum_{i=1}^l \frac{\|S_i\|_\alpha^2}{l^2} = O(1).$$

Since $\alpha/2 > 2$, $\{[A_l - \mathbb{E}(A_l)]^2, l \ge 1\}$ is uniformly integrable [Chow and Teicher (1988)]. Hence, the weak convergence of $A_l$ implies the $\mathcal{L}^2$ moment convergence

$$\mathbb{E}\{[A_l - \mathbb{E}(A_l)]^2\} \to \sigma^4 \mathbb{E}\left\{ \int_0^1 [\mathbb{B}(t)^2 - \mathbb{E}(\mathbb{B}(t)^2)] \, dt \right\}^2 = \frac{\sigma^4}{3}.$$

**A.4. Proof of Corollary 3.** Choose $d \in \mathbb{N}$ such that $2^{d-1} < N \le 2^d$. Using the same argument as in the proof of Theorem 2 [see (41)–(43) therein], we have for $1 \le a < b$ that

$$\|V_b - V_a - \mathbb{E}(V_b - V_a)\|_{\alpha/2} = \left\|\sum_{i=a+1}^{b}(W_i^2 - \mathbb{E}W_i^2)\right\|_{\alpha/2}$$
$$= O[b^{3(1-1/p)/2}(b-a)^{2/\alpha}],$$

where the constant in $O$ does not depend on $a$ and $b$. To show (21), we shall apply a useful maximal inequality established in Wu (2007). By Proposition 1 in the latter paper,

$$\left\|\max_{n \le 2^d}|V_n - \mathbb{E}V_n|\right\|_{\alpha/2}$$
$$\le \sum_{r=0}^{d}\left[\sum_{l=1}^{2^{d-r}}\|V_{2^r l} - V_{2^r(l-1)} - \mathbb{E}(V_{2^r l} - V_{2^r(l-1)})\|_{\alpha/2}^{\alpha/2}\right]^{2/\alpha}.$$

Note that

$$\sum_{l=1}^{2^{d-r}}\|V_{2^r l} - V_{2^r(l-1)} - \mathbb{E}(V_{2^r l} - V_{2^r(l-1)})\|_{\alpha/2}^{\alpha/2}$$
$$= \sum_{l=1}^{2^{d-r}}O\{[(2^r l)^{3(1-1/p)/2}(2^r)^{2/\alpha}]^{\alpha/2}\}$$
$$= O(1)2^{r+3r(1-1/p)\alpha/4}(2^{d-r})^{1+3(1-1/p)\alpha/4}.$$

Hence,

$$\left\|\max_{n \le 2^d}|V_n - \mathbb{E}V_n|\right\|_{\alpha/2} = O(d+1)(2^d)^{2/\alpha+3(1-1/p)/2}$$

and (21) follows in view of $2^{d-1} < N \le 2^d$.

We now show (22). Note that $\alpha/2 > 1$. By (21), we have

$$\sum_{d=1}^{\infty}\frac{\|\max_{n \le 2^d}|V_n - \mathbb{E}V_n|\|_{\alpha/2}^{\alpha/2}}{(2^{d\tau}d^2)^{\alpha/2}} = \sum_{d=1}^{\infty}O(d^{-\alpha/2}) < \infty,$$

which by the Borel–Cantelli lemma implies that $V_N - \mathbb{E}V_N = o[N^\tau(\log N)^2]$ almost surely. Consequently, (22) easily follows from $\mathbb{E}V_n - t_n\sigma^2 = O[n^{1+(1-q)(1-1/p)}]$.   □

**Acknowledgments.** I thank Peter Glynn and a referee for their many useful comments.

Department of Statistics
University of Chicago
5734 Street University Avenue
Chicago, Illinois 60637
USA
E-mail: wbwu@galton.uchicago.edu